\documentclass[12pt]{amsart} 
\usepackage{amsmath} \usepackage{amssymb} \usepackage{amsthm}
\usepackage{amscd} \usepackage[all]{xy} 
\def\a{{\mathfrak{a}}} \def\b{{\mathfrak{b}}} \def\C{{\mathfrak{c}}}
\def\F{{\mathbb{F}}} \def\m{{\mathfrak{m}}}\def\n{{\mathfrak{n}}}
\def\N{{\mathbb{N}}} 
\def\Q{{\mathbb{Q}}} \def\Hom{{\mathrm{Hom}}} \def\Ker{{\mathrm{Ker}}}
\def\Im{{\mathrm{Im}}}  \def\Ann{{\mathrm{Ann}}}
 
\theoremstyle{plain}
\newtheorem{thm}{Theorem}[section] \newtheorem{cor}[thm]{Corollary}
\newtheorem{prop}[thm]{Proposition}
\newtheorem{propdef}[thm]{Proposition-Definition} \newtheorem{lem}[thm]{Lemma}
\theoremstyle{definition} \newtheorem{dfn}[thm]{Definition}
\newtheorem{eg}[thm]{Example} \theoremstyle{remark}
\newtheorem{rmk}[thm]{Remark}

\title{On a generalization of test ideals}
\author{Nobuo Hara}
\address{Mathematical Institute, Tohoku University, Sendai 980-8578, Japan}
\email{hara@math.tohoku.ac.jp}
\author{Shunsuke Takagi}
\address{Graduate School of Mathematical Sciences, University of Tokyo,
3-8-1, Komaba, Meguro, Tokyo 153-8914, Japan}
\email{stakagi@ms.u-tokyo.ac.jp}
\thanks{Both authors thank MSRI for the support and hospitality during
their stay in the fall of 2002. The first-named author is partially 
supported by Grant-in-Aid for Scientific Research, Japan.}
 
\begin{document}
\begin{abstract}
The test ideal $\tau(R)$ of a ring $R$ of prime characteristic is an 
important object in the theory of tight closure. In this paper, we study 
a generalization of the test ideal, which is the ideal $\tau(\a^t)$ 
associated to a given ideal $\a$ with rational exponent $t \ge 0$. 
We first prove a key lemma of this paper (Lemma \ref{key lemma}), which 
gives a characterization of the ideal $\tau(\a^t)$. As applications of 
this key lemma, we generalize the preceding results on the  behavior of 
the test ideal $\tau(R)$. Moreover, we prove an analog of so-called 
Skoda's theorem, which is formulated algebraically via adjoint ideals 
by Lipman in his proof of the ``modified Brian\c con--Skoda theorem."
\end{abstract}
 
\maketitle
 
\section*{Introduction}
 
Let $R$ be a Noetherian commutative ring of characteristic $p > 0$. 
The test ideal $\tau(R)$ of $R$, introduced by Hochster and Huneke 
\cite{HH}, is defined to be the annihilator ideal of all tight closure 
relations in $R$ and plays an important role in the theory of tight 
closure. In \cite{HY}, the first-named author and Yoshida introduced 
a generalization of tight closure, which we call $\a$-tight closure 
associated to a given ideal $\a$, and defined the ideal $\tau(\a)$ 
to be the annihilator ideal of all $\a$-tight closure relations in $R$. 
We can also consider $\a^t$-tight closure and the ideal $\tau(\a^t)$ 
with rational exponent $t \ge 0$ (or, rational coefficient in a sense), 
and even more, those with several rational exponents. 

The ideals $\tau(\a^t)$ have several nice properties similar to those 
of multiplier ideals ${\mathcal J}(\a^t)$ defined via resolution of 
singularities in characteristic zero; see \cite{La} for a systematic 
study of multiplier ideals. Among them, we have an analog of Lipman's 
``modified Brian\c con--Skoda theorem" (\cite[Theorem 2.1]{HY}, cf.\ 
\cite{Li}) and the subadditivity theorem in regular local rings 
(\cite[Theorem 4.5]{HY}, cf.\ \cite{DEL}). It is notable that the above 
properties of the ideals $\tau(\a^t)$ are proved quite algebraically 
via characteristic $p$ methods. On the other hand, we can 
prove that the multiplier ideal ${\mathcal J}(\a^t)$ in a normal 
$\Q$-Gorenstein ring of characteristic zero coincides, after reduction 
to characteristic $p \gg 0$, with the ideal $\tau(\a^t)$; see \cite{HY} 
and also \cite{T}.

In this paper, we study further properties of $\a^t$-tight closure 
and the ideal $\tau(\a^t)$ in characteristic $p > 0$, improve results 
obtained in \cite{HY}, and give some applications. To do this, we first 
prove a key lemma of this paper (Lemma \ref{key lemma}), which gives a 
characterization of (elements of) the ideal $\tau(\a^t)$. 
Precisely speaking, Lemma \ref{key lemma} characterizes the ideal 
$\tilde{\tau}(\a^t)$ given in Definition \ref{tildetau}, which is 
contained in the ideal $\tau(\a^t)$ and expected to coincide with 
$\tau(\a^t)$. The identification of the ideals $\tau(\a^t)$ and 
$\tilde{\tau}(\a^t)$ holds true in some reasonable situations, 
e.g., in normal $\Q$-Gorenstein rings; see \cite{AM}, \cite{HY} and  \cite{LS}. 

Although the description of Lemma \ref{key lemma} is somewhat 
complicated, it turns out to be very useful in studying various 
properties of the ideals $\tau(\a^t)$. As applications of this key 
lemma, we answer a question raised in \cite{HY} (Corollary \ref{HYQ}) 
and consider the relationship of test elements and $\a^t$-test elements 
(Corollary \ref{test}). We also apply Lemma \ref{key lemma} to the 
study of the behavior of the ideals $\tau(\a^t)$ under localization 
(Proposition \ref{loc}), completion (Proposition \ref{com}) 
and finite morphisms which are \'etale in codimension one (Theorem 
\ref{etale}). 
These results generalize the preceding results \cite{LS} and \cite{BS} 
on the behavior of the test ideal $\tau(R)$, and we hope that the proofs 
become simpler with the use of Lemma \ref{key lemma}. 

Other ingredients of this paper are Theorems \ref{skoda1} and \ref{skoda2}, 
which assert that if $\a$ is an ideal with a reduction generated by $l$ 
elements, then $\tau(\a^l) = \tau(\a^{l-1})\a$. 
This is an analog of so-called Skoda's theorem \cite{La}, which is 
formulated algebraically via adjoint ideals (or multiplier ideals) 
by Lipman \cite{Li} in his proof of the ``modified Brian\c con--Skoda 
theorem." Theorem \ref{skoda1} is a refinement of \cite[Theorem 2.1]{HY} 
and is proved by an easy observation on the relationship of regular powers 
and Frobenius powers of ideals (cf.\ \cite{AH}). Theorem 4.2 is based on 
the same idea, but is proved under a slightly different assumption with 
the aid of Lemma \ref{key lemma}. We can find an advantage of the ideal 
$\tau(\a)$ in the simplicity of the proofs of Theorems \ref{skoda1} and 
\ref{skoda2}, because the proof of Skoda's theorem for multiplier ideals 
needs a deep vanishing theorem which is proved only in characteristic zero 
\cite{La}, \cite{Li}.

\section{Preliminaries}

In this paper, all rings are excellent reduced commutative rings with unity. 
For a ring $R$, we denote by $R^{\circ}$ the set of elements of $R$ which are 
not in any minimal prime ideal. Let $R$ be a ring of prime characteristic $p >
0$ and $F\colon R \to R$ the Frobenius map which sends $x \in R$ to $x^p \in 
R$. The ring $R$ viewed as an $R$-module via the $e$-times iterated Frobenius 
map $F^e \colon R \to R$ is denoted by ${}^e\! R$. Since $R$ is assumed to be 
reduced, we can identify $F \colon R \to {}^e\! R$ with the natural inclusion 
map $R \hookrightarrow R^{1/p^e}$. We say that $R$ is {\it F-finite} if ${}^1
\! R$ (or $R^{1/p}$) is a finitely generated $R$-module. 

Let $R$ be a ring of characteristic $p > 0$ and let $M$ be an $R$-module. For 
each $e \in \N$, we denote $\F^e(M) = \F_R^e(M) := {}^e\! R \otimes_R M$ and 
regard it as an $R$-module by the action of $R = {}^e\! R$ from the left. 
Then we have the induced $e$-times iterated Frobenius map $F^e \colon M \to 
\F^e(M)$. The image of $z \in M$ via this map is denoted by $z^q:= F^e(z) 
\in \F^e(M)$. For an $R$-submodule $N$ of $M$, we denote by $N^{[q]}_M$ the 
image of the induced map $\F^e(N) \to \F^e(M)$. 

Now we recall the definition of $\a^t$-tight closure. See \cite{HY} 
for details.

\begin{dfn}\label{ta}
Let $\a$ be an ideal of a ring $R$ of characteristic $p>0$ such that $\a \cap R^{\circ} \ne \emptyset$, and let 
$N \subseteq M$ be $R$-modules. Given a rational number $t \ge 0$, the 
{\it $\a^t$-tight closure} $N^{*\a^t}_M$ 
of $N$ in $M$ is defined to be the submodule of $M$ consisting of all 
elements $z \in M$ for which there exists $c \in R^{\circ}$ such that 
$$cz^q\a^{\lceil tq \rceil} \subseteq N^{[q]}_M$$ 
for all large $q = p^e$, where $\lceil tq \rceil$ is the least integer 
which is greater than or equal to $tq$. The $\a^t$-tight closure of an ideal $I \subseteq R$ is just 
defined by $I^{*\a^t} = I^{*\a^t}_R$. 
\end{dfn}

\begin{rmk}\label{a^t}
The rational exponent $t$ for $\a^t$-tight closure in Definition 1.1 is just 
a formal notation, but it is compatible with ``real" powers of the ideal. 
Namely, if $\b = \a^n$ for $n \in \N$, then $\a^t$-tight closure is the 
same as $\b^{t/n}$-tight closure. This allows us to extend the definition 
to several rational exponents: Given ideals $\a_1,\dots ,\a_r \subseteq R$ 
with $\a_i \cap R^{\circ} \ne \emptyset$ 
and rational numbers $t_1,\dots ,t_r \ge 0$, if $t_i = tn_i$ for nonnegative 
$t \in \Q$ and $n_i \in \N$ with $i = 1,\dots ,r$, we can define $\a_1^{t_1}
\cdots\a_r^{t_r}$-tight closure to be $(\a_1^{n_1}\cdots\a_r^{n_r})^t
$-tight closure. If $N$ is a submodule of an $R$-module $M$, then an 
element $z \in M$ is in the $\a_1^{t_1}\cdots\a_r^{t_r}$-tight closure 
$N^{*\a_1^{t_1}\cdots\a_r^{t_r}}_M$ of $N$ in $M$ if and only if there 
exists $c \in R^{\circ}$ such that $cz^q\a_1^{\lceil t_1q \rceil}\cdots
\a_r^{\lceil t_rq \rceil} \subseteq N^{[q]}_M$ for all $q = p^e \gg 0$.
\end{rmk}

Since $N^{*\a^t}_M/N \cong 0^{*\a^t}_{M/N}$ for $R$-modules $N \subseteq M$ 
(\cite[Proposition 1.3 (1)]{HY}), the case where $N = 0$ is essential. 
Using the $\a^t$-tight closure of the zero submodule, we can define two 
ideals $\tau(\a^t)$ and $\tilde{\tau}(\a^t)$.

\begin{propdef}\label{taudef}
{\rm (\cite[Definition-Theorem 6.5]{HY})}
Let $R$ be an excellent reduced ring of characteristic $p > 0$, $\a \subseteq R$ an ideal such that $\a \cap R^{\circ} \ne \emptyset$ and $t \ge 0$ a rational number. Let $E =\bigoplus_{\m}
E_R(R/\m)$ be the direct sum, taken over all maximal ideals $\m$ of $R$, 
of the injective envelopes of the residue field $R/\m$. Then the following
ideals are equal to each other and we denote it by $\tau(\a^t)$. 
\begin{enumerate}
\renewcommand{\labelenumi}{(\roman{enumi})}
\item $\displaystyle\bigcap_M \Ann_R(0^{*\a^t}_M)$, where $M$ runs 
through all finitely generated $R$-modules. 
\vspace{3pt}
\item $\displaystyle\bigcap_{M\subseteq E}\!\!\Ann_R(0^{*\a^t}_M)$, where 
$M$~runs through all finitely generated submodules of~$E$. 
\item $\displaystyle\bigcap_{J\subseteq R} (J:J^{*\a^t})$, where $J$ 
runs through all ideals of $R$. 
\end{enumerate}
\end{propdef}

The description (ii) of $\tau(\a^t)$ in Proposition-Definition \ref{taudef} means that 
$\tau(\a^t)$ is the annihilator ideal of the ``finitistic $\a^t$-tight 
closure" of zero in $E$, that is, the union of $0^{*\a^t}_M$ in $E$ 
taken over all finitely generated submodules $M$ of $E$. It would be 
natural to consider the $\a^t$-tight closure of zero in $E$ instead 
of the finitistic tight closure. 

\begin{dfn}\label{tildetau}
Let $E =\underset{\m}{\bigoplus} E_R(R/\m)$ be as in Proposition-Definition $\ref{taudef}$. 
Then we define the ideal $\tilde{\tau}(\a^t)$ by 
$$\tilde{\tau}(\a^t) = \Ann_R(0_E^{*\a^t}).$$
We have $\tilde{\tau}(\a^t) \subseteq \tau(\a^t)$ in general, since 
$0^{*\a^t}_M \subseteq 0^{*\a^t}_E$ for all (finitely generated) submodules 
of $E$. We do not know any example in which the $\a^t$-tight closure of 
zero and the finitistic $\a^t$-tight closure of zero disagree in $E$, 
and it seems reasonable to assume the following condition. 

\vspace{8pt}\noindent
(1.4.1) \quad
We say that {\it condition $(*)$ is satisfied for $\a^t$} if 
$\tau(\a^t) = \tilde{\tau}(\a^t)$. 
\end{dfn}

Condition $(*)$ is satisfied in many situations. For example, if 
$R$ is a graded ring, then condition $(*)$ is satisfied for the 
unit ideal $R$ (see \cite{LS1}), and if $R$ is an excellent $\Q$-Gorenstein normal local ring, then condition $(*)$ is satisfied for every rational 
number $t \ge 0$ and every ideal $\a \subseteq R$ such that $\a \cap R^{\circ} \ne \emptyset$ (\cite[Theorem 1.13]
{HY}; see also \cite{AM}). 

In most part of this paper, we will assume condition $(*)$ to deduce 
various properties of the ideal $\tau(\a^t)$, because what we actually do is to 
prove those properties for the ideal $\tilde{\tau}(\a^t)$, which are 
translated into those for $\tau(\a^t)$ under condition $(*)$. So, we 
remark that most of our results do hold true for $\tilde{\tau}(\a^t)$ 
without assuming condition $(*)$. 

Finally we recall the notion of an $\a^t$-test element, which is quite 
useful to treat the $\a^t$-tight closure operation. See \cite[Theorems 
1.7 and 6.4]{HY} for the existence of $\a^t$-test elements. 

\begin{dfn}\label{testdef}
Let $\a$ be an ideal of a ring $R$ of characteristic $p > 0$ such that $\a \cap R^{\circ} \ne \emptyset$, and let 
$t \ge 0$ be a rational number. An element $d \in R^{\circ}$ is called 
an {\it $\a^t$-test element} if for every finitely generated $R$-module 
$M$ and $z \in M$, the following holds: $z \in 0_M^{*\a^t}$ if and only 
if $dz^q\a^{\lceil tq \rceil} = 0$ for all powers $q = p^e$ of $p$.
\end{dfn}

\begin{rmk}\label{testrem} 
In the case where $\a = R$ is the unit ideal, the ideal $\tau(\a) 
= \tau(R)$ is called the test ideal of $R$ and an $R$-test element 
is nothing but a test element as defined in \cite{HH}. In this case, 
$\tau(R) \cap R^{\circ}$ is exactly equal to the set of test elements. 
However, $\tau(\a^t) \cap R^{\circ}$ is not equal to the set of 
$\a^t$-test elements in general. The relationship between the 
ideal $\tau(\a^t)$ and $\a^t$-test elements is not a priori clear, 
but we will see later in Corollary \ref{test} that an element of 
$\tau(\a^t) \cap R^{\circ}$ is always an $\a^t$-test element. 
%
\end{rmk} 

The reader is referred to \cite{HY} for basic properties of $\a^t$-tight 
closure and the ideal $\tau(\a^t)$. Among them, see especially Propositions 
1.3 and 1.11 and Theorem 1.7, together with Section 6, of \cite{HY}.  
%

All results in this paper are proved in characteristic $p > 0$. But 
it may help better understanding of the results to keep in mind the 
correspondence of the ideals $\tau(\a^t)$ and the multiplier ideals 
defined via resolution of singularities in characteristic zero 
(\cite{HY}). Namely, given a rational number $t \ge 0$ and an ideal 
$\a$ of a normal $\Q$-Gorenstein ring essentially of finite type over 
a field of characteristic zero, the multiplier ideal ${\mathcal J}(\a^t)$ 
coincides, after reduction to characteristic $p \gg 0$, with the ideal 
$\tau(a^t)$.

\section{A characterization of the ideal $\tau(\a^t)$}

We first show a key lemma of this paper, which gives a characterization 
of the ideal $\tilde{\tau}(\a^t)$ and reveals the relationship between the ideal $\tau(\a^t)$ and $\a^t$-test elements under condition $(*)$.
Its essential idea is found in 
\cite[Theorem 3.13]{T}.

\begin{lem}\label{key lemma}
Let $(R,\m)$ be an F-finite local ring of characteristic $p > 0$, 
$\a \subseteq R$ an ideal such that $\a \cap R^{\circ} \ne \emptyset$ and $t \ge 0$ a rational number.
Fix a system of generators $x^{(e)}_1,\ldots, x^{(e)}_{r_{e}}$ of $\a^{\lceil tq \rceil}$ for each $q = p^e$. Then for 
any element $c \in R$, the following three conditions are equivalent to 
each other.
\begin{enumerate}
\renewcommand{\labelenumi}{(\roman{enumi})}
\item $c \in \tilde{\tau}(\a^t)$.
\item
For any element $d \in R^{\circ}$ and any integer $e_0 \ge 0$, 
there exist an integer $e_1 \geq e_0$ and $R$-homomorphisms 
$\phi^{(e)}_i \in \Hom_R(R^{1/p^e},R)$ for $e_0 \le e \le e_1$ and 
$1 \le i \le r_e$ such that 
$$c=\sum_{e=e_0}^{e_1}\sum_{i=1}^{r_e} \phi^{(e)}_i((dx_i^{(e)})^{1/p^e}).$$
\item
For an $\a^t$-test element $d \in R^{\circ}$, there exist a positive 
integer $e_1$ and $R$-linear maps $\phi^{(e)}_i \in \Hom_R(R^{1/p^e},R)$ 
for $0 \le e \le e_1$ and $1 \le i \le r_e$ such that 
$$c=\sum_{e=0}^{e_1}\sum_{i=1}^{r_e} \phi^{(e)}_i((dx_i^{(e)})^{1/p^e}).$$
\item[(iii)']
For an $\a^t$-test element $d \in R^{\circ}$, there exist a positive 
integer $e_1$ and $R$-linear maps $\phi^{(e)}_i \in \Hom_R(R^{1/p^e},R)$ 
for $0 \le e \le e_1$ and $1 \le i \le r_e$ such that 
$$c \in \sum_{e=0}^{e_1}\sum_{i=1}^{r_e} \phi^{(e)}_i((d\a^{\lceil tp^e \rceil})^{1/p^e}).$$
\end{enumerate}
\end{lem}

\begin{proof}
We prove the equivalence of conditions (i) and (ii). 
Let $F^e \colon E_R \to {}^e\! R \otimes_R E_R$ be the $e$-times iterated 
Frobenius map induced on the injective envelope $E_R := E_R(R/\m)$ of the 
residue field $R/\m$. For an element $d \in R^{\circ}$ and an integer $e \ge 0$, 
we write 
$$d\mathbf{x}^{(e)}F^e := {}^t(dx^{(e)}_1F^e,\dots,dx^{(e)}_{r_e}F^e) 
                       \colon E_R \to ({}^e\! R \otimes_R E_R)^{\oplus r_e}.$$
By definition and condition $(*)$, $c \in \tau(\a^t)$ if and only if 
$c\cdot\bigcap_{e\ge e_0} \Ker (d\mathbf{x}^{(e)}F^e) = 0$ holds for 
every element $d \in R^{\circ}$ and every integer $e_0 \ge 0$. Since 
$E_R$ is an Artinian $R$-module, there exists an integer $e_1 \ge e_0$ 
(depending on $d \in R^{\circ}$ and $e_0 \ge 0$) such that 
$$\bigcap_{e\ge e_0} \Ker(d\mathbf{x}^{(e)}F^e) 
        = \bigcap_{e=e_0}^{e_1} \Ker(d\mathbf{x}^{(e)}F^e).$$
Therefore denoting
$$\Phi = \Phi_d^{(e_0,e_1)} := 
  {}^t(d\mathbf{x}^{(e_0)}F^{e_0},\ldots,d\mathbf{x}^{(e_1)}F^{e_1})
        \colon E_R \to 
        \bigoplus_{e=e_0}^{e_1} ({}^e\! R \otimes_R E_R)^{\oplus r_e},$$ 
we see that $c \in \tau(\a^t)$ if and only if for every $d \in R^{\circ}$ and 
$e_0\ge 0$, there exists $e_1\ge e_0$ such that $c\cdot\Ker(\Phi_d^{(e_0,e_1)})
= 0$.

Since $R$ is F-finite, the map $\Phi = \Phi_d^{(e_0,e_1)} \colon E_R \to \bigoplus_{e=e_0}^{e_1}
({}^e\! R \otimes_R E_R)^{\oplus r_e}$ is the Matlis dual of the map 
$$\Psi = (\psi^{(e_0)},\ldots,\psi^{(e_1)}) \colon 
          \bigoplus_{e=e_0}^{e_1} \Hom_R(R^{1/p^e},R)^{\oplus r_e} \to R,$$ 
where the map $\psi^{(e)} \colon \Hom_R(R^{1/p^e},R)^{\oplus r_e} \to \Hom_R
(R,R) = R$ is induced by the $R$-linear map $R \to (R^{1/p^e})^{\oplus r_e}$ 
sending $1$ to $((dx_1^{(e)})^{1/p^e},\ldots,(dx_{r_e}^{(e)})^{1/p^e})$. 
It then follows that $c \cdot \Ker(\Phi) = 0$ if and only if $c \in \Im(\Psi)$.
By the definition of $\Psi$, this is equivalent to saying that there exist 
$R$-linear maps $\phi^{(e)}_i \in \mathrm{Hom}_R(R^{1/p^e},R)$ for $e_0 \le 
e \le e_1$ and $1 \le i \le r_e$ such that
$$c=\sum_{e=e_0}^{e_1}\sum_{i=1}^{r_e} \phi^{(e)}_i((dx_i^{(e)})^{1/p^e}).$$

The equivalence of conditions (i) and (iii) (resp. (iii)') is obtained in the same way. 
\end{proof}

\begin{rmk}\label{mixed key lemma}
An advantage of Lemma \ref{key lemma} is that it is applicable even in the 
absense of the completeness of $R$. For example, compare the hypotheses in 
Theorem \ref{skoda1} and Theorem \ref{skoda2}.
%
\end{rmk}

As immediate consequences of Lemma \ref{key lemma}, we have the following 
corollaries. First Lemma \ref{key lemma} gives an affirmative answer to a 
question raised in \cite[Discussion 5.18]{HY}.

\begin{cor}\label{HYQ}
Let $(R,\m)$ be an F-finite local ring of characteristic $p>0$, and 
suppose that condition $(*)$ is satisfied for the maximal ideal $\m$. 
Then $\tau(\m) = R$ if and only if for every $\m$-primary ideal 
$\a \subset R$, we have a strict containment $\tau(\a) \supsetneq \a$. 
In particular, if $R$ is a regular local ring with $\dim R \ge 2$, then $\tau(\a) \supsetneq \a$ for every 
$\m$-primary ideal $\a \subset R$.
\end{cor}

\begin{proof}
Suppose that $\tau(\m) = R$. We will show that $\tau(\a) \supseteq 
(\a:\m)$ for every ideal $\a \subseteq R$ such that $\a \cap R^{\circ} \ne \emptyset$. By Lemma \ref{key lemma}, 
for any element $d \in R^{\circ}$ and any integer $e_0 \ge 0$, there 
exist an integer $e\ge e_0$, an element $x \in \m^{p^e}$ and an 
$R$-module homomorphism $\phi \colon R^{1/p^e} \to R$ such that 
$\phi((dx)^{1/p^e})=1$. Now fix any element $y \in (\a:\m)$. 
Since $\phi((dxy^{p^e})^{1/p^e}) = \phi((dx)^{1/p^e})y = y$ and 
$xy^{p^e} \in \a^{p^e}$, we have $y \in \tau(\a)$ by using Lemma 
\ref{key lemma} again. The converse implication is trivial. The 
latter assertion follows from \cite[Theorem 2.15]{HY}.
\end{proof}

\begin{cor}\label{test}
Let $(R,\m)$ be an F-finite local ring of characteristic $p > 0$. 
If condition $(*)$ is satisfied for the unit ideal $R$, then a test 
element of $R$ is an $\a^t$-test element for all ideals $\a \subseteq R$ such that $\a \cap R^{\circ} \ne \emptyset$ 
and all rational numbers $t \ge 0$. 
\end{cor}

\begin{proof}
Let $c$ be a test element of $R$, that is, an element of $\tau(R) \cap 
R^{\circ}$. Given any ideals $\a,J$ of $R$, any element $z \in J^{*\a^t}$ 
and any power $q$ of $p$, it is enough to show that 
$cz^q\a^{\lceil tq \rceil} \subseteq J^{[q]}$. 
Since $z \in J^{*\a^t}$, there exist $d \in R^{\circ}$ and $e_0 \in \N$ 
such that $dz^Q\a^{\lceil tQ \rceil} \subseteq J^{[Q]}$ for every power 
$Q \ge p^{e_0}$ of $p$. Then by Lemma \ref{key lemma}, there exist 
$e_1 \in \N$ and $\phi_{e} \in \Hom_R(R^{1/p^e},R)$ for $e_0 \le e 
\le e_1$ such that $c = \sum_{e=e_0}^{e_1} \phi_e(d^{1/p^e})$. Since 
$dz^{qp^e}(\a^{\lceil tq \rceil})^{[p^e]} \subseteq 
        dz^{qp^e}\a^{\lceil tqp^e \rceil} \subseteq J^{[qp^e]}$ 
for every $e_0 \le e \le e_1$, we have 
$$d^{1/p^e}z^q\a^{\lceil tq \rceil}R^{1/p^e} \subseteq J^{[q]}R^{1/p^e}.$$
Applying $\phi_e$ and summing up, we obtain 
$$cz^q\a^{\lceil tq \rceil} = \sum_{e=e_0}^{e_1} \phi_e(d^{1/p^e})
z^q\a^{\lceil tq \rceil} \subseteq 
\sum_{e=e_0}^{e_1} \phi_e(J^{[q]}R^{1/p^e}) \subseteq J^{[q]},$$
as required. 
\end{proof}

\begin{rmk}\label{testexist}
Since $\tilde{\tau}(\a^t) \subseteq \tilde{\tau}(R)$, we see from Corollary \ref{test} 
that any element of $\tau(\a^t) \cap R^{\circ}$ is an $\a^t$-test element 
as long as condition $(*)$ is satisfied for $\a^t$. Also, Corollary 
\ref{test} is considered a refinement of \cite[Theorem 1.7]{HY}, 
which asserts that, if the localized ring $R_c$ at an element $c \in 
R^{\circ}$ is strongly F-regular, then some power $c^n$ of $c$ is an 
$\a^t$-test element for all ideals $\a \subseteq R$ such that $\a \cap R^{\circ} \ne \emptyset$ and all rational 
numbers $t \ge 0$. Indeed, if $R_c$ is strongly F-regular, then some 
power $c^n$ of $c$ is a test element (precisely speaking, an element of $\tilde{\tau}(R) \cap R^{\circ}$) by \cite{HH1} (see also \cite{HH2}), 
so that $c^n$ is an $\a^t$-test element for all $\a$ and $t$ by Corollary 
\ref{test}. 
%
\end{rmk}

\section{Behavior of $\tau(\a)$ under localization, completion and finite homomorphisms}

\begin{prop}\label{loc} {\rm (cf.\ \cite{LS})}
Let $(R,\m)$ be an F-finite local ring of characteristic $p >0$, $\a \subseteq R$ an ideal such that $\a \cap R^{\circ} \ne \emptyset$ and $t \ge 0$ a rational number. 
Let $W$ be a multiplicatively closed subset of $R$ and suppose that condition $(*)$ is satisfied for $\a^t$ and $(\a R_W)^t$. Then
$$\tau((\a R_W)^t) = \tau(\a^t)R_W.$$
\end{prop}

\begin{proof}
Fix a system of generators $x^{(e)}_1, \ldots, x^{(e)}_{r_e}$ of 
$\a^{\lceil tq \rceil}$ for each $q=p^e$. If an element $c \in R_W$ 
is contained in $\tau((\a R_W)^t)$, then by Lemma \ref{key lemma}, 
for any element $d \in R^{\circ}$ and any nonnegative integer $e_0$, 
there exist an integer $e_1 \ge e_0$ and $R_W$-homomorphisms 
$\phi^{(e)}_i \in \Hom_{R_W}(R_W^{1/p^e},R_W)$ for $e_0 \le e \le e_1$ 
and $1 \le i \le r_e$ such that 
$$c=\sum_{e=e_0}^{e_1}\sum_{i=1}^{r_e} \phi^{(e)}_i((dx_i^{(e)})^{1/p^e}).$$
Since $R$ is F-finite, there exists an element $y \in W$ such that we 
can regard $y\phi^{(e)}_i$ as an element of $\Hom_R(R^{1/p^e},R)$ for 
all $e_0 \le e \le e_1$ and $1 \le i \le r_e$. Therefore, thanks to 
Lemma \ref{key lemma} again, we have $cy \in \tau(\a^t)$, that is, 
$c \in \tau(\a^t)R_W$.
The converse argument just reverses this. The proposition is proved.
\end{proof}

\begin{prop}\label{com} {\rm (cf.\ \cite{LS})}
Let $(R,\m)$ be an F-finite local ring of characteristic $p>0$, 
$\a \subseteq R$ an ideal such that $\a \cap R^{\circ} \ne \emptyset$ and $t \ge 0$ a rational number. 
Let $\widehat{R}$ denote the $\m$-adic completion of $R$ and suppose that condition $(*)$ is satisfied for $\a^t$ and $(a\widehat{R})^t$. Then
$$\tau((\a\widehat{R})^t)=\tau(\a^t)\widehat{R}.$$
\end{prop}

\begin{proof}
Fix a system of generators $x^{(e)}_1, \ldots, x^{(e)}_{r_e}$ of 
$\a^{\lceil tq \rceil}$ for each $q=p^e$. 
Since $R$ is F-finite, we can take an element $d \in R^{\circ}$ which is an $\a^t$- and $(\a\widehat{R})^t$-test element by \cite[Theorem 6.4]{HY}.
If an element $c \in \widehat{R}$ 
is contained in $\tau((\a \widehat{R})^t)$, then by Lemma \ref{key lemma},  
there exist an integer $e_1 > 0$ and $\widehat{R}$-homomorphisms 
$\phi^{(e)}_i \in \Hom_{\widehat{R}}(\widehat{R}^{1/p^e},\widehat{R})$ for $0 \le e \le e_1$ 
and $1 \le i \le r_e$ such that 
$$c=\sum_{e=0}^{e_1}\sum_{i=1}^{r_e} \phi^{(e)}_i((dx_i^{(e)})^{1/p^e}).$$
Since $R$ is F-finite, 
$\Hom_{\widehat{R}}(\widehat{R}^{1/p^e},\widehat{R}) \cong \widehat{R} \otimes_R \Hom_R(R^{1/p^e}, R)$, so that 
there exist $y^{(e)}_{i,1}, \dots, y^{(e)}_{i,s_{e,i}} \in \widehat{R}$ and $\psi^{(e)}_{i,1}, \dots, \psi^{(e)}_{i,s_{e,i}} \in \Hom_R(R^{1/p^e},R)$ with 
$\phi^{(e)}_i= \sum_{j=1}^{s_{e,i}} y^{(e)}_{i,j} \otimes \psi^{(e)}_{i,j}$ for  all $0 \le e \le e_1$ and $1 \le i \le r_e$. 
Then
$$c=\sum_{e=0}^{e_1}\sum_{i=1}^{r_e}\sum_{j=1}^{s_{e,i}} y^{(e)}_{i,j} \psi^{(e)}_{i,j}((dx_i^{(e)})^{1/p^e}).$$
Therefore, thanks to 
Lemma \ref{key lemma} again, we have $c \in \tau(\a^t)\widehat{R}$.
The converse argument just reverses this. The proposition is proved.
\end{proof}

Our characterization is also applicable to the following situation; 
see also \cite{BS}.

\begin{thm}\label{etale}
Let $(R,\m) \hookrightarrow (S,\n)$ be a pure finite local homomorphism 
of F-finite normal local rings of characteristic $p >0$ which is \'{e}tale 
in codimension one. Let $\a$ be a nonzero ideal of $R$ and let $t$ be a nonnegative 
rational number. Assume that condition $(*)$ is satisfied for $\a^t$ and 
$(\a S)^t$. Then 
$$\tau((\a S)^t) \cap R=\tau(\a^t).$$
Moreover if $R \hookrightarrow S$ is flat, then 
$$\tau((\a S)^t)=\tau(\a^t)S.$$
\end{thm}

\begin{proof}
Since $R \hookrightarrow S$ is pure, we have $\tau((\a S)^t) \cap R 
\subseteq \tau(\a^t)$ by \cite[Proposition 1.12]{HY}. We will prove 
the reverse inclusion. Fix a system of generators $x^{(e)}_1, \ldots, 
x^{(e)}_{r_e}$ of $\a^{\lceil tq \rceil}$ for every $q=p^e$. 
Take an element $d \in R^{\circ}$ which is a common test element for $R$ and $S$ (cf. \cite[Remark 6.5]{BS}). 
Then $d$ is not only an $\a^t$-test element but also an $(\a S)^t$-test element by Corollary \ref{test} and Remark \ref{testexist}.
If an element $c$ belongs to $\tau(\a^t)$, then by Lemma \ref{key lemma}, 
there exist an integer $e_1 > 0$ and $R$-linear maps 
$\phi^{(e)}_i\in\Hom_R(R^{1/p^e},R)$ for $0 \le e \le e_1$ and 
$1 \le i \le r_e$ such that
$$c=\sum_{e=0}^{e_1}\sum_{i=1}^{r_e} \phi^{(e)}_i((dx_i^{(e)})^{1/p^e}).$$
Let $c_{e,i}=\phi^{(e)}_i((dx_i^{(e)})^{1/p^e})$ and let $E_S=E_S(S/\n)$ be 
the injective envelope of the residue field of $S$. Since $R \hookrightarrow 
S$ is \'{e}tale in codimension one, by tensoring $\phi^{(e)}_i$ with $E_S$ 
over $R$, we have the commutative diagram 
$$\xymatrix
{E_S \ar[r]^{dx_i^{(e)}F^e_S \hspace{0.7cm}} \ar[dr]_{c_{e,i}} & {}^e\! S \otimes_S E_S \ar[d] \\
& E_S, \\
}$$
where $F^e_S \colon E_S \to {}^e\! S \otimes_S E_S$ is the induced 
$e$-times iterated Frobenius map on $E_S$ (see \cite[Proof of Theorem 2.7]{W}). 
Hence $\Ker(dx_i^{(e)}F^e_S) \subseteq (0:c_{e,i})_{E_S}$ for every 
$0 \le e \le {e_1}$ and $1 \le i \le r_e$. Since 
$\sum_{e=0}^{e_1}\sum_{i=1}^{r_e} c_{e,i} = c$, we have 
$c \cdot \bigcap_{e=0}^{e_1} \bigcap_{i=1}^{r_e} \Ker(dx_i^{(e)}F^e_S) = 0$.
By condition $(*)$, this implies that $c \in \tau((\a S)^t)$. 

Now we will show the latter assertion assuming that $R \hookrightarrow S$ 
is flat. By the above argument, we already know that $\tau((\a S)^t) 
\supseteq \tau(\a^t)S$. Suppose that $c \in \tau((\a S)^t)$. Since $S$ is 
a free $R$-module, we can choose a basis $s_1, \ldots, s_k$ for $S$ over $R$ and write $c = \sum_{j=1}^k c_js_j$ for $c_j \in R$. 
For any ideal $I \subseteq R$ and any $z \in I^{*\a^t}$, clearly 
$z \in (IS)^{*(\a S)^t}$. 
By the definition of $\tau((\a S)^t)$, 
we have $cz = \sum_{j=1}^k c_j z s_j \in IS = \bigoplus_{j=1}^k I s_j$. 
It follows that $c_jz \in I$, therefore $c_j \in \tau(\a^t)$ for every 
$j=1, \ldots, k$. Thus $c \in \tau(\a^t)S$.
\end{proof}

The following example shows that the last equality in Theorem \ref{etale}
breaks down in the absence of the flatness, even for the case $\a = R$.

\begin{eg}
Let $S = k[[x,y,z]]/(x^n+y^n+z^n)$ be the Fermat hypersurface of degree $n$ over a field of characteristic $p \ge n+1$, 
$R = S^{(r)}$ the $r$th Veronese subring of $S$, 
and assume that $n \ge 3$ and $r \ge 2$ are not divisible by $p$. 
Then 
$\tau(R) = R_{\ge\lfloor 1+(n-3)/r\rfloor} = R_{\ge\lceil (n-2)/r\rceil}$
(note that $\lfloor 1+(n-3)/r \rfloor = \lfloor (n-2)/r+(r-1)/r \rfloor 
= \lceil (n-2)/r \rceil$) and $\tau(S) = S_{\ge n-2}$. Hence $\tau(S) 
\cap R = R_{\ge \lceil (n-2)/r \rceil} = \tau(R)$, but $\tau(R)S = 
S_{\ge \lceil (n-2)/r \rceil r} \subsetneq S_{\ge n-2} = \tau(S)$ if 
$n-2$ is not divisible by $r$.
\end{eg}

\section{Lipman-Skoda's theorem}

In \cite{Li}, Lipman proves under the Grauert--Riemenschneider vanishing
theorem that, if $\a$ is an ideal of a regular local ring with a reduction 
generated by $l$ elements, then ${\mathcal J}(\a^l)={\mathcal J}(\a^{l-1})
\a$, where ${\mathcal J}(\b)$ denotes the multiplier ideal (or adjoint ideal 
in the sense of \cite{Li}) associated to an ideal $\b$. This result is 
called Skoda's theorem \cite{La} and formulated algebraically by Lipman 
in his proof of ``modified Brian\c con--Skoda theorem." We give a simple 
proof of the corresponding equality for the ideal $\tau(\a^l)$. 
The following version is just a refinement of \cite[Theorem 2.1]{HY}; 
see also Remark \ref{a^t} for the definition of the ideal $\tau(\a^l\b^t)$ 
with ``bi-exponents."

\begin{thm}\label{skoda1}
Let $(R,\m)$ be a complete local ring of characteristic $p > 0$ and 
let $\a \subseteq R$ be an ideal of positive height with a reduction 
generated by $l$ elements. Let $\b$ be an ideal of $R$ such that $\b \cap R^{\circ} \ne \emptyset$ and $t \ge 0$ a rational number. Then
$$\tau(\a^l \b^t) = \tau(\a^{l-1} \b^t)\a.$$
\end{thm}

\begin{proof}
We will see that
$$0^{*\a^l\b^t}_M = 0^{*\a^{l-1}\b^t}_M \colon \a \; {\rm in}\; M$$
for any $R$-module $M$. By fundamental properties of $\a^t$-tight
closure \cite[Proposition 1.3]{HY}, the inclusion $0^{*\a^l\b^t}_M
\subseteq 0^{*\a^{l-1}\b^t}_M \colon \a$ is immediate, and to prove 
the reverse inclusion, we may assume without loss of generality that 
$\a$ itself is generated by $l$ elements. Let $z \in 0^{*\a^{l-1}\b^t}_M 
\colon \a$, i.e., $z\a \subseteq 0^{*\a^{l-1}\b^t}_M$. Then there exists
$c \in R^{\circ}$ such that $cz^q\a^{[q]}(\a^{l-1})^q\b^{\lceil tq \rceil} 
= 0$ in $\F^e(M)$ for all $q = p^e \gg 0$. Since $\a$ is generated by 
$l$ elements, one has $\a^{ql} = \a^{q(l-1)}\a^{[q]}$, so that $cz^q
\a^{lq}\b^{\lceil tq \rceil} = 0$ for all $q = p^e \gg 0$, that is, 
$z \in 0^{*\a^l\b^t}_M$. Thus we have $0^{*\a^l\b^t}_M = 0^{*\a^{l-1}
\b^t}_M \colon \a$.

Now assume that $(R,\m)$ is a complete local ring and let $E = E_R(R/\m)$, 
the injective envelope of the $R$-module $R/\m$. Then by the Matlis 
duality, $\Ann_E(\tau(\a^{l-1}\b^t))$ is equal to the union of 
$0^{*\a^{l-1}\b^t}_M$ taken over all finitely generated $R$-submodules 
$M$ of $E$. Hence, if $z \in \Ann_E (\tau(\a^{l-1}\b^t)\a)$, then there 
exists a finitely generated submodule $M \subset E$ such that $z \in 
(0^{*\a^{l-1}\b^t}_M \colon \a)_E$. Replacing $M$ by $M+Rz \subset E$,
one has $z \in (0^{*\a^{l-1}\b^t}_M \colon \a)_M = 0^{*\a^l\b^t}_M$. 
Consequently, $\Ann_E (\tau(\a^{l-1}\b^t)\a)$ is equal to the union 
of $(0^{*\a^{l-1}\b^t}_M \colon \a)_M = 0^{*\a^l\b^t}_M$ taken over 
all finitely generated submodules $M \subset E$. Therefore
$$\tau(\a^l\b^t) 
            = \bigcap_{M\subset E} \Ann_R(0^{*\a^l\b^t}_M)
            = \Ann_R(\Ann_E(\tau(\a^{l-1}\b^t)\a)) = \tau(\a^{l-1}\b^t)\a.$$ 
\end{proof}

The characterization of the ideal $\tau(\a^t)$ given in Lemma \ref{key lemma} 
enables us to replace the completeness assumption in the above theorem by 
condition $(*)$.

\begin{thm}\label{skoda2}
Let $(R,\m)$ be an F-finite local ring of characteristic $p > 0 $ 
and let $\a \subseteq R$ be an ideal of positive height with a reduction generated by $l$ 
elements. 
Let $\b$ be an ideal of $R$ such that $\b \cap R^{\circ} \ne \emptyset$ and $t$ a nonnegative rational number. 
If condition $(*)$ is satisfied for $\a^{l-1}\b^t$ and $\a^l\b^t$, 
then
$$\tau(\a^l\b^t) = \tau(\a^{l-1}\b^t)\a.$$
\end{thm}

\begin{proof}
We may assume without loss of generality that $\a$ is generated by $l$
elements, so that $\a^{ql} = \a^{q(l-1)}\a^{[q]}$ for all $q = p^e$.
We fix an element $d \in R^{\circ}$ such that $R_d$ is regular. Then, 
thanks to Remark \ref{testexist}, some power $d^n$ is an $\a^k\b^t$-test 
element for all $k \ge 0$.

By Lemma \ref{key lemma}, 
an element $c \in R$ is in $\tau(\a^l\b^t)$ if and only if there exist 
finitely many $R$-homomorphisms $\phi^{(e)}_i \in\Hom_R(R^{1/p^e},R)$ 
for $0 \le e\le e_1$ and $1\le i\le r_e$ such that
$$c \in \sum_{e=0}^{e_1}\sum_{i=1}^{r_e}
        \phi^{(e)}_i((d^n\a^{p^el}\b^{\lceil tp^e \rceil})^{1/p^e}).$$
Since
\begin{align*}
\phi^{(e)}_i((d^n\a^{p^el}\b^{\lceil tp^e \rceil})^{1/p^e}) & 
= \phi^{(e)}_i((d^n\a^{p^e(l-1)}\a^{[p^e]}\b^{\lceil tp^e \rceil})^{1/p^e}) \\
 & = \phi^{(e)}_i((d^n\a^{p^e(l-1)}\b^{\lceil tp^e \rceil})^{1/p^e})\a \\
 & \subseteq \tau(\a^{l-1}\b^t)\a
\end{align*}
again by Lemma \ref{key lemma}, this is equivalent to saying that
$c \in \tau(\a^{l-1}\b^t)\a$.
\end{proof}

\begin{cor} {\rm (Modified Brian\c con--Skoda, cf.\ \cite{BSk}, \cite{HY}, 
\cite{HH}, \cite{Li})}
Let $(R,\m)$ and $\a \subseteq R$ be as in Theorem $\ref{skoda1}$ or 
$\ref{skoda2}$. Then 
$$\tau(\a^{n+l-1}) \subseteq \a^n$$ 
for all $n \ge 0$. In particular, if $R$ is weakly F-regular, then 
$\overline{\a^{n+l-1}} \subseteq \a^n$ for all $n \ge 0$. 
\end{cor}

\begin{cor}
Let $(R,\m)$ be a $d$-dimensional local ring of characteristic $p > 0$ 
with infinite residue field $R/\m$ and $\a \subseteq R$ an ideal of 
positive height. Let $\b$ be an ideal of $R$ such that $\b \cap R^{\circ} \ne \emptyset$ and $t$ a nonnegative rational number. 
We assume that $(R,\m)$ is complete or it is F-finite 
and condition $(*)$ is satisfied for $\a^{n+d-1}\b^t$ for every $n \ge 0$. 
Then for any $n \ge 0$ one has
$$\tau(\a^{n+d-1}\b^t) = \tau(\a^{d-1}\b^t)\a^n.$$
\end{cor}

\begin{proof}
We can assume that $\a$ is a proper ideal of $R$. Since the residue field 
$R/\m$ is infinite, by \cite{NR}, $\a$ has a reduction ideal generated by 
at most $d$ elements. Therefore the assertion immediately follows from 
Theorems \ref{skoda1} and \ref{skoda2}.
\end{proof}

\begin{eg}{\rm (cf.\ \cite[Theorem 2.15]{HY})}
Let $R$ be a $d$-dimensional regular local ring of characteristic $p>0$ with the maximal ideal $\m$. Then
$$\tau(\m^n)=\left\{ 
\begin{array}{ll}
R & \mbox{if $n < d$,} \\
\m^{n-d+1} & \mbox{if $n \ge d$.}
\end{array}
\right.$$
In particular, $\tau(\m^{d-1}) \supsetneq \tau(\m^{d-2})\m$. 
This shows that $l \ge d$ is the best possible bound for the equality $\tau(\m^l)=\tau(\m^{l-1})\m$ in this case.
\end{eg}

\begin{rmk}
Let $(R,\m)$ and $\a$ be as in Theorem \ref{skoda1} or \ref{skoda2}
and let $\mathfrak q$ be any reduction of $\a$. Then $\tau(\a^{l-1})\a 
= \tau(\a^{l-1}){\mathfrak q}$. In particular, $\tau(\a^{l-1}) = \tau
({\mathfrak q}^{l-1})$ is contained in the coefficient ideal $\C (\a,
{\mathfrak q})$ of $\a$ with respect to $\mathfrak q$; cf.\ \cite{AH}. 

If $R$ is Gorenstein, $\a$ is $\m$-primary with minimal reduction 
$\mathfrak q$ and if the Rees algebra $R[\a t]$ is F-rational, then 
the equality $\tau(\a^{d-1}) = \C(\a,{\mathfrak q})$ holds; see 
\cite{HY} and Hyry \cite{Hy}.
\end{rmk}

\end{document}